\input amstex
\input amsppt.sty
\magnification=\magstep1
\hsize=30truecc
\vsize=22.2truecm
\baselineskip=16truept
\TagsOnRight
\pageno=1
\nologo
\def\Z{\Bbb Z}
\def\N{\Bbb N}

\def\l{\left}
\def\r{\right}
\def\bg{\bigg}
\def\({\bg(}
\def\[{\bg\lfloor}
\def\){\bg)}
\def\]{\bg\rfloor}
\def\t{\text}
\def\f{\frac}

\def\B{B_{k+1}}

\def\bi{\binom}
\def\eq{\equiv}

\def\ls{\leqslant}
\def\gs{\geqslant}
\def\mo{\roman{mod}}

\def\ve{\varepsilon}
\def\al{\alpha}

\def\Proof{\noindent{\it Proof}}

\def\Remark{\medskip\noindent{\it  Remark}}

\hbox {Adv. in Appl. Math. 45(2010), no.\,1, 125--148.}
\bigskip
\topmatter
\title New congruences for central binomial coefficients\endtitle
\author Zhi-Wei Sun$^1$ and Roberto Tauraso$^2$\endauthor
\leftheadtext{Zhi-Wei Sun and Roberto Tauraso}
\affil $^1$Department of Mathematics, Nanjing University\\
 Nanjing 210093, People's Republic of China\\zwsun\@nju.edu.cn
\\{\tt http://math.nju.edu.cn/$\sim$zwsun}
\medskip
$^2$Dipartimento di Matematica
\\Universit\`a di Roma ``Tor Vergata"
\\Roma 00133, Italy
\\tauraso\@mat.uniroma2.it
\\{\tt http://www.mat.uniroma2.it/$\sim$tauraso}
\endaffil
\abstract Let $p$ be a prime and let $a$ be a positive integer. In this paper we determine
$\sum_{k=0}^{p^a-1}\bi{2k}{k+d}/m^k$ and
$\sum_{k=1}^{p-1}\bi{2k}{k+d}/(km^{k-1})$ modulo $p$ for all $d=0,\ldots,p^a$, where $m$ is any
integer not divisible by $p$. For example, we show that if
$p\not=2,5$ then
$$\sum_{k=1}^{p-1}(-1)^k\f{\bi{2k}k}k\eq-5\f{F_{p-(\f p5)}}p\ (\mo\ p),$$
where $F_n$ is the $n$th Fibonacci number and $(-)$ is the Jacobi symbol. We also prove that
if $p>3$ then
$$\sum_{k=1}^{p-1}\f{\bi{2k}k}k\eq\f 89p^2B_{p-3}\ (\mo\ p^3),$$
where $B_n$ denotes the $n$th Bernoulli number.
\endabstract
\thanks 2010 {\it Mathematics Subject Classification}.\,Primary 11B65;
Secondary 05A10, 05A19, 11A07, 11B39, 11B68.
\newline\indent The first author is the corresponding author, and he is
supported by the National Natural Science Foundation (grant 10871087)
and the Overseas Cooperation Fund (grant 10928101) of China.
\endthanks
\endtopmatter
\document

\heading{1. Introduction}\endheading

A central binomial coefficient has the form
$\bi{2n}n$ with $n\in\N=\{0,1,\ldots\}$.
A well-known theorem of Wolstenholme (see, e.g., [5]) states that
$$\f12\bi{2p}p=\bi{2p-1}{p-1}\eq1\ (\mo\ p^3)\quad\t{for any prime}\ p>3.$$
In 2006 H. Pan and Z. W. Sun [9]
used a sophisticated combinatorial identity to deduce that if $p$ is a prime then
$$\sum_{k=0}^{p-1}\bi{2k}{k+d}\eq\l(\f{p-d}3\r)\ (\mo\ p)\quad \t{for}\ d=0,\ldots,p,\tag1.1$$
where the Jacobi symbol $(\f a3)$ coincides with the unique integer $\ve\in\{0,\pm1\}$ satisfying $a\eq\ve\ (\mo\ 3)$.
In a recent paper [16] the authors determined $\sum_{k=0}^{p^a-1}\bi{2k}{k+d}$ mod $p^2$ for any prime $p$ and
$d\in\{0,1,\ldots,p^a\}$ with $a\in\Z^+=\{1,2,3,\ldots\}$.

 In this paper we extend the congruence (1.1) in a new way and derive various congruences related to recurrences.
 Throughout this paper, for an assertion $A$ we set
 $$[A]=\cases1&\t{if}\ A \ \t{holds},\\0&\t{otherwise}.\endcases$$
We also define two recurrences
$\{u_n(x)\}_{n\in\N}$ and  $\{v_n(x)\}_{n\in\N}$ of polynomials as follows:
$$u_0(x)=0,\ u_1(x)=1,\ \t{and}\ u_{n+1}(x)=xu_{n}(x)-u_{n-1}(x)\ (n=1,2,\ldots),$$
and
$$v_0(x)=2,\ v_1(x)=x,\ \t{and}\ v_{n+1}(x)=xv_n(x)-v_{n-1}(x)\ (n=1,2,\ldots).$$
For a fixed integer $x$, the sequences $\{u_n(x)\}_{n\in\N}$ and $\{v_n(x)\}_{n\in\N}$
are linear recurrences of integers.
By induction, for any $n\in\N$ we have
$$u_n(-x)=(-1)^{n-1}u_n(x)\quad\t{and}\quad v_n(-x)=(-1)^nv_n(x).\tag1.2$$

Now we state our first theorem.

\proclaim{Theorem 1.1} Let $p$ be a prime and let $d\in\{0,\ldots,p^a\}$ with $a\in\Z^+$.
 Let $m\in\Z$ with $p\nmid m$. Then we have
$$\sum_{k=0}^{p^a-1}\f{\bi{2k}{k+d}}{m^k}\eq u_{p^a-d}(m-2)\ (\mo\ p)\tag1.3$$
and $$d\sum_{k=1}^{p^a-1}\f{\bi{2k}{k+d}}{km^{k-1}}
\eq 2(-1)^{d}+v_{p^a-d}(m-2)\ (\mo\ p)\ \ \t{provided}\ d>0.\tag1.4$$
If $p\not=2$, then
$$\sum_{k=0}^{p^a-1}\f{\bi{2k}{k+d}}{m^k}\eq-u_{d-(\f{m(m-4)}{p^a})}(m-2)
\ (\mo\ p)\tag1.5$$
and also
$$d\sum_{k=1}^{p^a-1}\f{\bi{2k}{k+d}}{km^{k-1}}
\eq 2(-1)^{d}+v_{d-(\f{m(m-4)}{p^a})}(m-2)\ (\mo\ p)
\ \t{provided}\ d>0,
\tag1.6$$
where $u_{-1}(x)=xu_0(x)-u_1(x)=-1$ and
$v_{-1}(x)=xv_0(x)-v_1(x)=x$.
\endproclaim

\Remark\ 1.1. Let $p$ be any prime and let $a\in\Z^+$. As $u_n(-1)=(\f n3)$ for
$n=0,1,2,\ldots$,  (1.3) in the case $m=1$ yields that
$$\sum_{k=0}^{p^a-1}\bi{2k}{k+d}\eq\l(\f{p^a-d}3\r)\ (\mo\ p)\quad \t{for every}\ d=0,1,\ldots,p^a.$$
Since $v_n(-1)=3[3\mid n]-1$ for all $n\in\N$, by (1.4) in the case $m=1$, for $d\in\{1,\ldots,p^a\}$ we have
$$d\sum_{k=1}^{p^a-1}\f{\bi{2k}{k+d}}k\eq
\cases 2(-1)^d+2\ (\mo\ p)&\t{if}\ p^a\eq d\ (\mo\ 3),\\2(-1)^d-1\ (\mo\ p)&\t{otherwise}.
\endcases$$
\bigskip

The well-known Fibonacci sequence $\{F_n\}_{n\in\N}$ is defined by
$$F_0=0,\ F_1=1,\ \t{and}\ F_{n+1}=F_n+F_{n-1}\ \t{for}\ n=1,2,3,\ldots.$$
Its companion $\{L_n\}_{n\in\N}$, the Lucas sequence, is given by
$$L_0=2,\ L_1=1,\ \t{and}\ L_{n+1}=L_n+L_{n-1}\ \t{for}\ n=1,2,3,\ldots.$$
Define
$$\align&F_{-1}=F_1-F_0=1,\ F_{-2}=F_0-F_{-1}=-1,
\\&L_{-1}=L_1-L_0=-1,\ L_{-2}=L_0-L_{-1}=3.\endalign$$
By induction,
$F_{2n}=u_n(3)$ and $L_{2n}=v_n(3)$ for $n=-1,0,1,\ldots$.
Note also that $u_{2n}(0)=v_{2n+1}(0)=0$ and $v_{2n}(0)/2=u_{2n+1}(0)=(-1)^n$ for all $n\in\N$.
Thus, with the help of (1.2), Theorem 1.1 in the cases $m=-1,2$
gives the following consequence.

\proclaim{Corollary 1.1} Let $p$ be an odd prime and let
$d\in\{0,1,\ldots,p^a\}$ with $a\in\Z^+$. Then
$$\sum_{k=0}^{p^a-1}(-1)^k\bi{2k}{k+d}\eq(-1)^{d-[p\not=5]}F_{2(d-(\f{p^a}5))}\ (\mo\ p),\tag1.7$$
and
$$d\sum_{k=1}^{p^a-1}(-1)^k\f{\bi{2k}{k+d}}k
\eq(-1)^{d-[p=5]}L_{2(d-(\f{p^a}5))}-2(-1)^d\ (\mo\ p)\tag1.8$$
provided $d>0$. Also,
$$\sum_{k=0}^{p^a-1}\f{\bi{2k}{k+d}}{2^k}
\eq\cases0\ (\mo\ p)&\t{if}\ p^a\eq d\ (\mo\ 2),
\\1\ (\mo\ p)&\t{if}\ p^a\eq d+1\ (\mo\ 4),\\-1\ (\mo\ p)&\t{if}\ p^a\eq d-1\ (\mo\ 4),\endcases
\tag1.9$$
and for $d>0$ we have
$$d\sum_{k=1}^{p^a-1}\f{\bi{2k}{k+d}}{k2^{k}}-(-1)^d
\eq\cases0\ (\mo\ p)&\t{if}\ p^a\not\eq d\ (\mo\ 2),
\\1\ (\mo\ p)&\t{if}\ p^a\eq d\ (\mo\ 4),\\-1\ (\mo\ p)&\t{if}\ p^a\eq d+2\ (\mo\ 4).\endcases
\tag1.10$$
\endproclaim

Our following result can be viewed as a complement to Theorem 1.1.

\proclaim{Theorem 1.2} Let $p$ be a prime and let $m$ be an
integer not divisible by $p$. Then we have
$$\f12\sum_{k=1}^{p-1}(-1)^k\f{\bi{2k}k}{km^{k-1}}\eq\f{m^p-V_p(m)}p\
(\mo\ p),\tag1.11$$
where the polynomial sequence $\{V_n(x)\}_{n\in\N}$ is
defined as follows:
$$V_0(x)=2,\ V_1(x)=x,\ \t{and}\ V_{n+1}(x)=x(V_n(x)+V_{n-1}(x))\
(n\in\Z^+).$$
\endproclaim

Given a prime $p$ and an integer $a$ not divisible by $p$, we use $q_p(a)$ to denote the integer
$(a^{p-1}-1)/p$ and call $q_p(a)$ a {\it Fermat quotient} with base $a$.
See E. Lehmer [7] for connections between Fermat quotients and Fermat's last theorem.

\proclaim{Corollary 1.2} Let $p$ be an odd prime. Then
$$\sum_{k=1}^{p-1}\f{\bi {2k}k}{k2^{k-1}}\eq\sum_{k=1}^{p-1}\f{\bi{2k}k}{k4^{k}}\eq 2q_p(2)\ (\mo\ p).\tag1.12$$
If $p\not=3$ then
$$\sum_{k=1}^{p-1}\f{\bi{2k}k}{k3^{k-1}}\eq3q_p(3)\ (\mo\ p).\tag1.13$$
\endproclaim

\proclaim{Corollary 1.3} Let $p$ be an odd prime.

{\rm (i)} If $p\not=5$, then we have
$$\align\sum_{k=1}^{p-1}(-1)^k\f{\bi{2k}k}{k}\eq&-5\f{F_{p-(\f p5)}}p\ (\mo\ p),\tag1.14
\\\sum_{k=1}^{p-1}(-1)^k\f{\bi{2k}k}{k5^k}\eq&q_p(5)-6\f{F_{p-(\f p5)}}p\ (\mo\ p),\tag1.15
\\\sum_{k=1}^{p-1}\f{\bi{2k}k}{k5^k}\eq&q_p(5)-\f{F_{p-(\f p5)}}p\ (\mo\ p).\tag1.16
\endalign$$

{\rm (ii)} Define the Pell sequence $\{P_n\}_{n\in\N}$ by
$$P_0=0,\ P_1=1,\ \t{and}\ P_{n+1}=2P_n+P_{n-1}\ (n=1,2,3,\ldots).$$
Then
$$\sum_{k=1}^{p-1}(-1)^k\f{\bi{2k}k}{k4^k}\eq 2q_p(2)-4\f{P_{p-(\f2p)}}p
\eq2\sum_{0<k<3p/4}\f{(-1)^{k-1}}k\ (\mo\ p).\tag1.17$$

{\rm (iii)} Let $\{S_n\}_{n\in\N}$ be the sequence defined by
$$S_0=0,\ S_1=1,\ \t{and}\ S_{n+1}=4S_n-S_{n-1}\ (n=1,2,3,\ldots).$$
If $p>3$, then
$$\sum_{k=1}^{p-1}(-1)^k\f{\bi{2k}k}{k2^{k}}
\eq q_p(2)-6\l(\f 2p\r)\f{S_{(p-(\f 3p))/2}}p
\eq\sum_{0<k<5p/6}\f{(-1)^{k-1}}k\ (\mo\ p)\tag1.18$$
and
$$\sum_{k=1}^{p-1}\f{\bi{2k}k}{k6^k}\eq q_p(2)+q_p(3)-2\l(\f 2p\r)\f{S_{(p-(\f 3p))/2}}p\ (\mo\ p).
\tag 1.19$$
\endproclaim

\Remark\ 1.2. (a) A prime $p\not=2,5$ is called a Wall-Sun-Sun prime
if $F_{p-(\f p5)}\eq0\ (\mo\ p^2)$ (cf. [1]). In 1992 Z. H. Sun and Z. W. Sun [13]
showed that Fermat's equation $x^p+y^p=z^p$ has no integer solutions satisfying $p\nmid xyz$
unless $p$ is a Wall-Sun-Sun prime. There are no Wall-Sun-Sun primes
below $2\times 10^{14}$ (cf. [8]). In 1982 H. C. Williams [10] showed that
$$\f{F_{p-(\f p5)}}p\eq\f25\sum_{0<k<4p/5}\f{(-1)^k}k\ (\mo\ p).$$
(b) The second congruences in (1.17) and (1.18) are essentially due to
Z. W. Sun [14, 15]. For other information about the sequence $\{S_n\}_{n\in\N}$ the reader may consult [11].
\medskip

 In 2006 Pan and Sun [9] proved that
$$\sum_{k=1}^{p-1}\f{\bi{2k}k}k\eq0\ (\mo\ p)$$
for any prime $p>3$.
Here we determine the sum modulo $p^3$.

\proclaim{Theorem 1.3} Let $p$ be any prime and let $a\in\Z^+$. Then
we have
$$p^{a-1}\sum_{k=1}^{p^a-1}\f{\bi{2k}{k}}k\eq
\cases  2\ (\mo\ p^3)&\t{if}\ p=2,\\5\ (\mo\ p^3)&\t{if}\ p=3,
\\\f{8}{9}p^2B_{p-3}\ (\mo\ p^3)&\t{otherwise},\endcases\tag1.20$$
where $B_0,B_1,B_2,\ldots$ are the well-known Bernoulli numbers.
\endproclaim

 The following conjecture, which is related to (1.7) in the case $d=0$, seems very challenging.

 \proclaim{Conjecture 1.1} Let $p\not=2,5$ be a prime and let $a\in\Z^+$. Then
 $$\sum_{k=0}^{p^a-1}(-1)^k\bi{2k}k
\eq\l(\f {p^a}5\r)\l(1-2F_{p^a-(\f {p^a}5)}\r)\ (\mo\ p^3).$$
 \endproclaim

In the next section we are going to present two auxiliary identities.
Theorem 1.1, Theorem 1.2 and Corollaries 1.2-1.3,
and Theorem 1.3 will be proved in Sections 3, 4 and 5 respectively.

\heading{2. An auxiliary theorem}\endheading

\proclaim{Theorem 2.1} For any $n\in\Z^+$ and $d\in\Z$, we have
$$\aligned&\sum_{0\ls k<n}\bi{2k}{k+d}x^{n-1-k}+[d>0]x^nu_d(x-2)
\\&\ \ =\sum_{0\ls k<n+d}\bi{2n}ku_{n+d-k}(x-2)
\endaligned\tag2.1$$
and
$$\aligned&d\sum_{0<k<n}\f{\bi{2k}{k+d}}kx^{n-k}-[d\gs0]x^nv_{d}(x-2)
+[d=0]x^n
\\&\ =-\sum_{0\ls k<n+d}\bi{2n}kv_{n+d-k}(x-2)-2\bi{2n-1}{n+d-1}.
\endaligned\tag2.2$$
\endproclaim
\Proof. (i) We use induction on $n\in\Z^+$ to prove (2.1).

Since $(x-2)u_d(x-2)=u_{d+1}(x-2)+u_{d-1}(x-2)$ for
$d=1,2,3,\ldots$, we can easily see that (2.1) with $n=1$ holds for
all $d\in\Z$.

Now fix $n\in\Z^+$ and assume (2.1) for all $d\in\Z$. Let $d$ be any
integer. For $k\in\N$, it is easy to see that
$$\bi{2n+2}k=\bi{2n}k+2\bi{2n}{k-1}+\bi{2n}{k-2}.$$
Thus,
$$\align&\sum_{0\ls k<(n+1)+d}\bi{2n+2}ku_{n+1+d-k}(x-2)
\\=&\sum_{0\ls k<n+(d+1)}\bi{2n}ku_{n+(d+1)-k}(x-2)
\\&+2\sum_{0\ls j<n+d}\bi{2n}ju_{n+d-j}(x-2)
\\&+\sum_{0\ls i<n+(d-1)}\bi{2n}iu_{n+(d-1)-i}(x-2).
\endalign$$

By the induction hypothesis, for any $r\in\Z$ we have
$$\sum_{0\ls k<n+r}\bi{2n}ku_{n+r-k}(x-2)
=\sum_{0\ls k<n}\bi{2k}{k+r}x^{n-1-k}+[r>0]x^nu_r(x-2).$$
So, from the above we get
$$\align&\sum_{0\ls k<(n+1)+d}\bi{2n+2}ku_{n+1+d-k}(x-2)
\\=&\sum_{0\ls k<n}\(\bi{2k}{k+d+1}
+2\bi{2k}{k+d}+\bi{2k}{k+d-1}\)x^{n-1-k}
\\&+[d\gs0]x^nu_{d+1}(x-2)+2[d\gs0]x^nu_d(x-2)+[d>0]x^nu_{d-1}(x-2)
\\=&\sum_{0\ls k<n}\(\bi{2k+1}{k+d+1}+\bi{2k+1}{k+d}\)x^{n-1-k}-[d=0]x^nu_{-1}(x-2)
\\&+[d\gs0]x^n\l(u_{d+1}(x-2)+2u_d(x-2)+u_{d-1}(x-2)\r)
\\=&\sum_{0\ls k<n}\bi{2(k+1)}{(k+1)+d}x^{n-1-k}+[d=0]x^n+[d\gs0]x^nxu_d(x-2)
\\=&\sum_{0\ls k<n+1}\bi{2k}{k+d}x^{(n+1)-1-k}+[d>0]x^{n+1}u_d(x-2).
\endalign$$
This concludes the induction step and hence (2.1) holds.

(ii) By induction, $v_k(x-2)=2u_{k+1}(x-2)-(x-2)u_k(x-2)$ for all $k\in\Z$.
Thus, with the help of (2.1), we have
$$\align&\sum_{0\ls k\ls n+d}\bi{2n}kv_{n+d-k}(x-2)
\\=&2\sum_{0\ls k<n+d+1}\bi{2n}ku_{n+d+1-k}(x-2)
\\&-(x-2)\sum_{0\ls k<n+d}\bi{2n}ku_{n+d-k}(x-2)
\\=&2\sum_{0\ls k<n}\bi{2k}{k+d+1}x^{n-1-k}+[d+1>0]x^n2u_{d+1}(x-2)
\\&-(x-2)\(\sum_{0\ls k<n}\bi{2k}{k+d}x^{n-1-k}+[d>0]x^nu_d(x-2)\)
\\=&\sum_{0\ls k<n}\(2\bi{2k}{k+d+1}-(x-2)\bi{2k}{k+d}\)x^{n-1-k}+[d\gs0]x^nv_{d}(x-2).
\endalign$$

For $k\in\Z^+$ we have
$$\align&\bi{2k-2}{k+d}+\bi{2k-2}{k+d-1}=\bi{2k-1}{k+d}=\bi{2k-1}{k-d-1}
\\=&\f{k-d}{2k}\bi{2k}{k-d}=\f{k-d}{2k}\bi{2k}{k+d}=\f12\bi{2k}{k+d}-\f d{2k}\bi{2k}{k+d}.
\endalign$$
Thus
$$\align &\f12\sum_{0<k<n}\bi{2k}{k+d}x^{n-k}-\f d2\sum_{0<k<n}\f{\bi{2k}{k+d}}{k}x^{n-k}
\\=&\sum_{0<k\ls n}\(\bi{2k-2}{k+d}+\bi{2k-2}{k+d-1}\)x^{n-k}-\bi{(2n-2)+1}{n+d}
\\=&\sum_{0\ls k<n}\(\bi{2k}{k+d+1}+\bi{2k}{k+d}\)x^{n-1-k}-\bi{2n-1}{n+d}.
\endalign$$
It follows that
$$\align&d\sum_{0<k<n}\f{\bi{2k}{k+d}}kx^{n-k}+[d=0]x^n-2\bi{2n-1}{n+d}
\\&=\sum_{0\ls k<n}\((x-2)\bi{2k}{k+d}-2\bi{2k}{k+d+1}\)x^{n-1-k}.
\endalign$$

Combining the above we obtain
$$\align&\sum_{0\ls k\ls n+d}\bi{2n}kv_{n+d-k}(x-2)-[d\gs0]x^nv_d(x-2)
\\&\ =-d\sum_{0<k<n}\f{\bi{2k}{k+d}}kx^{n-k}-[d=0]x^n+2\bi{2n-1}{n+d},
\endalign$$
from which (2.2) follows.
\qed

\proclaim{Corollary 2.1} Let $n\in\Z^+$ and $d\in\N$. Then
$$\sum_{0\ls k<n}\bi{2k}{k+d}+\l(\f d3\r)=\sum_{0\ls k<n+d}\bi{2n}k\l(\f{n+d-k}3\r),\tag2.3$$
$$\sum_{0\ls k<n}(-1)^{k+d}\bi{2k}{k+d}+F_{2d}=\sum_{0\ls k<n+d}(-1)^k\bi{2n}k F_{2(n+d-k)},\tag2.4$$
and
$$\aligned &d\sum_{0<k<n}\f{(-1)^{k+d}}k\bi{2k}{k+d}+\sum_{0\ls k<n+d}\bi{2n}k(-1)^kL_{2(n+d-k)}
\\&\quad=L_{2d}-(-1)^{n+d}2\bi{2n-1}{n+d-1}-[d=0].\endaligned\tag2.5$$
\endproclaim
\Proof. For $j\in\N$ we have
$u_j(-1)=\l(\f j3\r)$, $(-1)^{j-1}u_j(-3)=u_j(3)=F_{2j}$ and $(-1)^jv_j(-3)=v_j(3)=L_{2j}$.
 Thus, (2.1) in the case $x=1$ yields (2.3), and
(2.1) and (2.2) in the case $x=-1$ reduce to (2.4) and (2.5) respectively.
This concludes the proof.
\qed

\heading{3. Proof of Theorem 1.1}\endheading

Given $A,B\in\Z$  we define the Lucas sequence $u_n=u_n(A,B)\ (n\in\N)$
and its companion $v_n=v_n(A,B)\ (n\in\N)$
as follows:
$$u_0=0,\ u_1=1,\ \t{and}\ u_{n+1}=Au_n-Bu_{n-1}\ \t{for}\ n=1,2,3,\ldots,$$
and
$$v_0=2,\ v_1=A,\ \t{and}\ v_{n+1}=Av_n-Bv_{n-1}\ \t{for}\ n=1,2,3,\ldots.$$
It is well known that
$$u_n=\sum_{0\ls k<n}\al^k\beta^{n-1-k}\ \ \t{and}\ \ v_n=\al^n+\beta^n
\quad\t{for all}\ n\in\N,$$
where $\al$ and $\beta$ are the two roots of the equation $x^2-Ax+B=0$.
It follows that if $n\in\N$ and  $m\in\{n,n+1,\ldots\}$ then
$$Au_n+v_n=2u_{n+1}\ \ \t{and}\ \ u_mv_n-u_nv_m=2B^nu_{m-n}.$$

\proclaim{Lemma 3.1} Let $A,B\in\Z$ with $B\not=0$. Let $u_n=u_n(A,B)$ for $n\in\N$, and define
$u_{-1}=(u_1-Au_0)/(-B)=-1/B$.
Let $p$ be an odd prime,  and let $a\in\Z^+$ and $d\in\{0,1,\ldots,p^a\}$.
Then we have
$$B^du_{p^a-d}\eq-c(A,B)u_{d-(\f{\Delta}{p^a})}\ (\mo\ p),\tag3.1$$
where $\Delta=A^2-4B$ and
$$c(A,B)=\cases A/2&\t{if}\ p\mid \Delta,\\B&\t{if}\ (\f{\Delta}{p^a})=1,\\1&\t{if}\ (\f{\Delta}{p^a})=-1.\endcases$$
\endproclaim
\Proof. The two roots of the equation $x^2-Ax+B=0$
are algebraic integers $\al=(A+\sqrt{\Delta})/2$ and $\beta=(A-\sqrt{\Delta})/2$.
Since
$$\bi{p^a}k=\f{p^a}k\bi{p^a-1}{k-1}\eq0\ (\mo\ p)\ \t{for}\ k=1,\ldots,p^a-1,$$
we have
 $$v_{p^a}=\al^{p^a}+\beta^{p^a}\eq(\al+\beta)^{p^a}=A^{p^a}\eq A^{p^{a-1}}\eq\cdots\eq A\ (\mo\ p)$$
 with the help of Fermat's little theorem.
 If $\Delta\not=0$, then
 $$\align u_{p^a}=&\f{\al^{p^a}-\beta^{p^a}}{\al-\beta}=\f1{\sqrt{\Delta}}\l(\l(\f{A+\sqrt{\Delta}}2\r)^{p^a}
 -\l(\f{A-\sqrt{\Delta}}2\r)^{p^a}\r)
 \\=&\f1{2^{p^a}\sqrt{\Delta}}\sum\Sb k=0\\2\nmid k\endSb^{p^a}\bi{p^a}k
 A^{p^a-k}\l((\sqrt{\Delta})^k-(-\sqrt{\Delta})^k\r)
 \\=&\f1{2^{p^a-1}}\sum^{p^a}\Sb k=1\\2\nmid k\endSb\bi{p^a}kA^{p^a-k}\Delta^{(k-1)/2};
 \endalign$$
 if $\Delta=0$ then $\al=\beta=A/2$ and hence $u_{p^a}=p^a(A/2)^{p^a-1}$. So we always have
 $$u_{p^a}=\f1{2^{p^a-1}}\sum^{p^a}\Sb k=1\\2\nmid k\endSb\bi{p^a}kA^{p^a-k}\Delta^{(k-1)/2}.$$
Note that $2^{p^a-1}\eq1\ (\mo\ p)$ by Fermat's little theorem. Thus, by Euler's criterion,
 $$u_{p^a}\eq\bi{p^a}{p^a}\Delta^{(p^a-1)/2}
=(\Delta^{(p-1)/2})^{\sum_{k=0}^{a-1}p^k}\eq\l(\f{\Delta}p\r)^a=\l(\f{\Delta}{p^a}\r)\ (\mo\ p).$$

 Observe that
$$2B^du_{p^a-d}=u_{p^a}v_d-u_dv_{p^a}\eq\l(\f{\Delta}{p^a}\r)v_d-u_dA\ (\mo\ p).$$
When $p\mid \Delta$, this yields
$$B^du_{p^a-d}\eq-\f A2u_d\ (\mo\ p).$$
If $(\f{\Delta}{p^a})=1$, then
$$2B^du_{p^a-d}\eq v_d-Au_d=2(u_{d+1}-Au_d)=-2Bu_{d-1}\ (\mo\ p)$$
and hence $B^du_{p^a-d}\eq-Bu_{d-1}\ (\mo\ p)$.
If $(\f{\Delta}{p^a})=-1$, then
$$2B^du_{p^a-d}\eq -v_d-Au_d=-2u_{d+1}\ (\mo\ p)$$
and thus $B^du_{p^a-d}\eq-u_{d+1}\ (\mo\ p)$. So (3.1) follows. \qed

\medskip
\noindent{\it Proof of Theorem 1.1}.
For $n=-1,0,1,\ldots$ let $u_n=u_n(m-2)$ and $v_n=v_n(m-2)$.

By Theorem 2.1,
$$\sum_{k=0}^{p^a-1}\bi{2k}{k-d}m^{p^a-1-k}=\sum_{0\ls
k<p^a-d}\bi{2p^a}ku_{p^a-d-k};$$
also, for $d>0$ we have
$$-d\sum_{0<k<p^a}\f{\bi{2k}{k-d}}km^{p^a-k}
=-\sum_{0\ls k<p^a-d}\bi{2p^a}k v_{p^a-d-k}-2\bi{2p^a-1}{p^a-d-1}.$$
By Fermat's little theorem, $m^{p^a}\eq m\ (\mo\ p)$. For
$k\in\{1,\ldots,p^a-1\}$ clearly
$$\bi{2p^a}k=\f{2p^a}k\bi{2p^a-1}{k-1}\eq0\ (\mo\ p);$$
also, if $d<p^a$ then
$$\bi{2p^a-1}{p^a-d-1}=\prod_{0<j<p^a-d}\l(\f{2p^a}j-1\r)\eq(-1)^{p^a-d-1}\eq(-1)^d\ (\mo\ p).$$
Therefore
$$\sum_{k=0}^{p^a-1}\f{\bi{2k}{k+d}}{m^k}\eq[d\not=p^a]\bi{2p^a}0u_{p^a-d}=u_{p^a-d}\ (\mo\ p);$$
if $d>0$ then
$$\align d\sum_{k=1}^{p^a-1}\f{\bi{2k}{k+d}}{km^{k-1}}
\eq&[d\not=p^a]\bi{2p^a}0v_{p^a-d}+2[d\not=p^a](-1)^d
\\\eq& v_{p^a-d}+2(-1)^d\ (\mo\ p).
\endalign$$
So we have (1.3) and (1.4).

 Now assume $p\not=2$ and set $\Delta=(m-2)^2-4\times1=m(m-4)$.
As $p\nmid m$, if $p\mid \Delta$
then $m\eq4\ (\mo\ p)$ and hence
$(m-2)/2\eq1\ (\mo\ p)$.
Thus, with the help of Lemma 3.1, we have
 $$\sum_{k=0}^{p^a-1}\f{\bi{2k}{k+d}}{m^k}\eq u_{p^a-d}\eq-u_{d-(\f{\Delta}{p^a})}
 \ (\mo\ p),$$
which proves (1.5). If $d>0$, then
 $$\align v_{d-(\f{\Delta}{p^a})}&=2u_{d-(\f{\Delta}{p^a})+1}-(m-2)u_{d-(\f{\Delta}{p^a})}
 \\=&-2u_{d-1-(\f{\Delta}{p^a})}+(m-2)u_{d-(\f{\Delta}{p^a})}
 \\\eq&2u_{p^a-d+1}-(m-2)u_{p^a-d}=v_{p^a-d} \ (\mo\ p).
 \endalign$$
Thus (1.6) follows from (1.4). We are done. \qed

 \heading{4. Proofs of Theorem 1.2 and Corollaries 1.2-1.3}
 \endheading

 \proclaim{Lemma 4.1} For any positive integer $n$, we have
$$\f12\sum_{0<k<n}\f{\bi{2k}k}{kx^k}=\sum_{0<d<n}(-1)^{d-1}
\sum_{0<k<n}\f{\bi{2k}{k+d}}{kx^k}.\tag4.1$$
\endproclaim
\Proof. Observe that
$$\align&\sum_{d=0}^{n-1}(-1)^{d}\sum_{0<k<n}\f{\bi{2k}{k+d}}{kx^k}
\\=&\sum_{0<k<n}\f1{k(-x)^k}\sum_{d=0}^{n-1}(-1)^{k+d}\bi{2k}{k+d}
\\=&\sum_{0<k<n}\f1{2k(-x)^k}\sum_{j=k}^{2k}\((-1)^j\bi{2k}j+(-1)^{2k-j}\bi{2k}{2k-j}\)
\\=&\sum_{0<k<n}\f1{2k(-x)^k}\l((1-1)^{2k}+(-1)^k\bi{2k}k\r)
=\f12\sum_{0<k<n}\f{\bi{2k}k}{kx^k}.
\endalign$$
So (4.1) follows. \qed

 \medskip
\noindent{\it Proof of Theorem 1.2}. By Lemma 4.1,
$$\f12\sum_{k=1}^{p-1}(-1)^k\f{\bi{2k}k}{km^{k-1}}
=\sum_{d=1}^{p-1}(-1)^{d}\sum_{k=1}^{p-1}\f{\bi{2k}{k+d}}{k(-m)^{k-1}}.$$

 In view of (1.4) and the basic fact
$$\f1p\bi pd=\f1d\prod_{0<k<d}\f{p-k}k\eq\f{(-1)^{d-1}}d\ (\mo\ p)\ \ (d=1,\ldots,p-1),$$
we have
 $$\align&\sum_{d=1}^{p-1}(-1)^{d}\sum_{k=1}^{p-1}\f{\bi{2k}{k+d}}{k(-m)^{k-1}}
 \\\eq&\sum_{d=1}^{p-1}\f{(-1)^{d}}d(v_{p-d}(-m-2)+2(-1)^d)
 \\\eq&\sum_{d=1}^{p-1}\f{(-1)^{d}}dv_{p-d}(-m-2)+\sum_{d=1}^{p-1}\l(\f1d+\f1{p-d}\r)
 \\\eq&-\f1p\sum_{d=1}^{p-1}\bi
 pdv_{p-d}(-m-2)=-\f1p\sum_{k=1}^{p-1}\bi pkv_k(-m-2)\ (\mo\ p).
 \endalign$$

 Let $\al$ and $\beta$ be the two roots of the equation
 $x^2-mx-m=0$. Then $(-\al-1)+(-\beta-1)=-m-2$ and
 $(-\al-1)(-\beta-1)=1$, also
 $$V_p(m)=\al^p+\beta^p\eq(\al+\beta)^p=m^p\eq m\ (\mo\ p).$$
 In the case $p\not=2$, we have
 $$\align&\sum_{k=1}^{p-1}\bi pkv_k(-m-2)=\sum_{k=1}^{p-1}\bi
 pk((-\al-1)^k+(-\beta-1)^k)
 \\=&(-\al)^p+(-\beta)^p-2-(-\al-1)^p-(-\beta-1)^p
 \\=&(-1)^pV_p(m)-2-(-1)^p\f{\al^{2p}+\beta^{2p}}{m^p}
\\=&-V_p(m)+\f{(\al^p+\beta^p)^2}{m^p}
 =\l(1+\f{V_p(m)-m^p}{m^p}\r)(V_p(m)-m^p)
 \\\eq& V_p(m)-m^p\ (\mo\ p^2)\quad (\t{since}\ V_p(m)\eq m^p\ (\mo\ p)).
 \endalign$$
 Note also that
$$\sum_{k=1}^{2-1}\bi 2kv_k(-m-2)=2(-m-2)\eq2m=V_2(m)-m^2\ (\mo\
2^2).$$
 Therefore (1.11) follows from the above. \qed

\medskip
\noindent{\it Proof of Corollary 1.2}. By induction, whenever
$n\in\N$ we have
$$\align &V_{4n}(-2)=(-1)^n2^{2n+1},\ V_{4n+1}(-2)=(-1)^{n+1}2^{2n+1},
\\&V_{4n+2}(-2)=0,\ V_{4n+3}(-2)=(-1)^n2^{2n+2}.\endalign$$
It follows that
$$V_p(-2)=-\l(\f 2p\r)2^{(p+1)/2}.$$
Combining this with (1.11) in the case $m=-2$, we get
$$\align\sum_{k=1}^{p-1}\f{\bi{2k}k}{k2^k}\eq&\f{V_p(-2)-(-2)^p}p=2^{(p+1)/2}\f{2^{(p-1)/2}-(\f 2p)}p
\\\eq&\(2^{(p-1)/2}+\l(\f 2p\r)\)\f{2^{(p-1)/2}-(\f 2p)}p=q_p(2)\ (\mo\ p).
\endalign$$
By induction, $V_n(-4)=(-1)^n2^{n+1}$ for all $n\in\N$. Thus, by
(1.11) with $m=-4$, we have
$$\f12\sum_{k=1}^{p-1}\f{\bi{2k}k}{k4^{k-1}}\eq\f{V_p(-4)-(-4)^p}p=2^p\f{2^p-2}p\eq 4q_p(2)\ (\mo\ p).$$
Therefore (1.12) holds.

 Now assume that $p\not=3$. By induction, for $n\in\N$ we have
$$V_n(-3)=\cases(3[3\mid n]-1)(-3)^{n/2}&\t{if}\ 2\mid n,
\\(\f n3)(-3)^{(n+1)/2}&\t{if}\ 2\nmid n.\endcases$$
Applying (1.11) with $m=-3$ we get
$$\align\f12\sum_{k=1}^{p-1}\f{\bi{2k}k}{k3^{k-1}}\eq&\f{V_p(-3)-(-3)^p}p=-(-3)^{(p+1)/2}\f{(-3)^{(p-1)/2}-(\f {-3}p)}p
\\\eq&\f32\((-3)^{(p-1)/2}+\l(\f{-3}p\r)\)\f{(-3)^{(p-1)/2}-(\f {-3}p)}p
\\\eq&\f 32\cdot\f{(-3)^{p-1}-1}p=\f 32q_p(3)\ (\mo\ p).
\endalign$$
So (1.13) is valid. \qed

 \medskip
\noindent{\it Proof of Corollary 1.3}.
(i) Applying Theorem 1.2 with $m=1$, we obtain that
$$\f12\sum_{k=1}^{p-1}(-1)^k\f{\bi{2k}k}k\eq\f{1-L_p}p\ (\mo\ p).$$
Let $\al$ and $\beta$ be the two roots of the equation $x^2-x-1=0$.
Suppose $p\not=5$ and set $n=(p-(\f p5))/2$. It is known that
$$L_n^2-5F_n^2=(\al^n+\beta^n)^2-(\al-\beta)^2\l(\f{\al^n-\beta^n}{\al-\beta}\r)^2
=4(\al\beta)^n=4(-1)^n$$
and
$$L_{2n}=\al^{2n}+\beta^{2n}=(\al^n+\beta^n)^2-2(\al\beta)^n=L_n^2-2(-1)^n.$$
By [13, Corollary 1],  $p\mid F_n$ if $p\eq1\ (\mo\ 4)$, and $p\mid L_n$ if $p\eq3\ (\mo\ 4)$.
Thus
$$L_{p-(\f p5)}=L_{2n}=5F_n^2+2(-1)^n=L_n^2-2(-1)^n\eq2\l(\f p5\r)\ (\mo\ p^2).$$
By induction,
$$2L_k=5F_{k-1}+L_{k-1}=5F_{k+1}-L_{k+1}\ \t{for}\ k=1,2,3,\ldots.$$
Therefore
$$2L_p=5F_{p-(\f p5)}+\l(\f p5\r)L_{p-(\f p5)}\eq5F_{p-(\f p5)}+2\ (\mo\ p^2)$$
and hence
$$\sum_{k=1}^{p-1}(-1)^k\f{\bi{2k}k}k\eq-2\f{L_p-1}p\eq-5\f{F_{p-(\f p5)}}p\ (\mo\ p).$$
This proves (1.14).

By (1.11) in the case $m=5$,
$$\sum_{k=1}^{p-1}(-1)^k\f{\bi{2k}k}{k5^k}\eq\f25\cdot\f{5^p-V_p(5)}p\ (\mo\ p).$$
Since $(5+3\sqrt5)/2$ and $(5-3\sqrt5)/2$ are the two roots of the equation $x^2-5x-5=0$,
$$\align V_p(5)=&\l(\f{5+3\sqrt5}2\r)^p+\l(\f{5-3\sqrt5}2\r)^p
\\=&\sqrt 5^p\l(\l(\f{1+\sqrt5}2\r)^{2p}-\l(\f{1-\sqrt5}2\r)^{2p}\r)
\\=&5^{(p+1)/2}\f{\al^p-\beta^p}{\al-\beta}
(\al^p+\beta^p)=5^{(p+1)/2}F_pL_p.
\endalign$$
As
$$L_p\eq1+\f52F_{p-(\f 5p)}\ (\mo\ p^2)$$
and
$$L_p=F_p+2F_{p-1}=2F_{p+1}-F_p=2F_{p-(\f p5)}+\l(\f p5\r)F_p,$$
we have
$$\align&\l(\f p5\r)F_pL_p=L_p(L_p-2F_{p-(\f p5)})
\\\eq&\l(1+\f52F_{p-(\f p5)}\r)\l(1+\f12F_{p-(\f p5)}\r)\eq1+3F_{p-(\f p5)}\ (\mo\ p^2)
\endalign$$
and hence
$$\align &V_p(5)=5^{(p+1)/2}F_pL_p
\\\eq&5^{(p+1)/2}\l(\f 5p\r)(1+3F_{p-(\f p5)})
\eq5^{(p+1)/2}\l(\f 5p\r)+15F_{p-(\f p5)}\ (\mo\ p^2).\endalign$$
Therefore
$$\align\sum_{k=1}^{p-1}(-1)^k\f{\bi{2k}k}{k5^k}\eq&\f25\cdot\f{5^p-5^{(p+1)/2}(\f 5p)-15F_{p-(\f p5)}}p
\\\eq&\l(5^{(p-1)/2}+\l(\f 5p\r)\r)\f{5^{(p-1)/2}-(\f 5p)}p-6\f{F_{p-(\f p5)}}p
\\\eq&q_p(5)-6\f{F_{p-(\f p5)}}p\ (\mo\ p).
\endalign$$
So (1.15) also holds.

Applying (1.11) with $m=-5$ we get
$$\f12\sum_{k=1}^{p-1}\f{\bi{2k}k}{k5^{k-1}}\eq\f{V_p(-5)+5^p}p\ (\mo\ p).$$
As the two roots of the equation $x^2+5x+5=0$ are $(-5\pm\sqrt5)/2$, we have
$$\align V_p(-5)=&\(\f{-5+\sqrt5}2\)^p+\(\f{-5-\sqrt5}2\)^p
\\=&\sqrt5^p\(\(\f{1-\sqrt5}2\)^p-\(\f{1+\sqrt5}2\)^p\)=-\sqrt5^{p+1}F_p.
\endalign$$
Recall that
$$\l(\f 5p\r)F_p=L_p-2F_{p-(\f p5)}\eq1+\f12F_{p-(\f p5)}\ (\mo\ p^2).$$
Thus
$$\align 5^{(p-1)/2}F_p-1\eq& 5^{(p-1)/2}\l(\f 5p\r)\l(1+\f12F_{p-(\f p5)}\r)-1
\\\eq&\l(\f 5p\r)\l(5^{(p-1)/2}-\l(\f 5p\r)\r)+\f12 F_{p-(\f p5)}
\\\eq&\f12\l(5^{(p-1)/2}+\l(\f 5p\r)\r)\l(5^{(p-1)/2}-\l(\f 5p\r)\r)+\f12F_{p-(\f p5)}
\\\eq&\f{5^{p-1}-1}2+\f12 F_{p-(\f p5)}\ (\mo\ p^2)
\endalign$$
and hence
$$\align \f12\sum_{k=1}^{p-1}\f{\bi{2k}k}{k5^{k-1}}\eq&\f{5^p-5^{(p+1)/2}F_p}p=\f{5^p-5}p-5\f{5^{(p-1)/2}F_p-1}p
\\\eq&5\(q_p(5)-\f{q_p(5)}2-\f{F_{p-(\f p5)}}{2p}\)
\ (\mo\ p).\endalign$$
This proves (1.16).

(ii) As $2+2\sqrt2$ and $2-2\sqrt2$ are the two roots of the equation $x^2-4x-4=0$, we have
$$V_p(4)=(2+2\sqrt2)^p+(2-2\sqrt2)^p=2^p\l((1+\sqrt2)^p+(1-\sqrt2)^p\r)=2^pQ_p,$$
where the sequence $\{Q_n\}_{n\in\N}$ is given by
$$Q_0=Q_1=2\ \t{and}\ Q_{n+1}=2Q_n+Q_{n-1}\ (n=1,2,3,\ldots).$$
By [15, Remark 3.1],
$$4\l(\f2p\r)P_p-Q_p=\l(\f 2p\r)Q_{p-(\f 2p)}\eq2\ (\mo\ p^2)$$
and $$P_{p-(\f 2p)}\eq\l(\f 2p\r)P_p-1\ (\mo\ p^2).$$
Thus
$$Q_p-2\eq4\l(\l(\f 2p\r)P_p-1\r)\eq4P_{p-(\f 2p)}\ (\mo\ p^2)$$
and hence
$$\align\sum_{k=1}^{p-1}(-1)^k\f{\bi{2k}k}{k4^k}\eq&\f{4^p-V_p(4)}{2p}
=2^{p-1}\f{2^p-Q_p}p
\\\eq&2q_p(2)-\f{Q_p-2}p\eq 2q_p(2)-4\f{P_{p-(\f 2p)}}p\ (\mo\ p)
\endalign$$
with the help of (1.11) in the case $m=4$.

By [14],
$$-2^{(p+1)/2}\f{P_p-2^{(p-1)/2}}p\eq\sum_{k=1}^{(p-1)/2}\f1{k2^k}
\eq\sum_{k=1}^{\lfloor3p/4\rfloor}\f{(-1)^{k-1}}k\ (\mo\ p).$$
(The last congruence was first conjectured by Z. H. Sun in 1988.)
Observe that
$$\align-2^{(p+1)/2}\f{P_p-2^{(p-1)/2}}p
\eq&-2^{(p+1)/2}\f{(\f2p)(1+P_{p-(\f2p)})-2^{(p-1)/2}}p
\\\eq&-2\f{P_{p-(\f 2p)}}p+2^{(p+1)/2}\f{2^{(p-1)/2}-(\f2p)}p
\\\eq&-2\f{P_{p-(\f 2p)}}p+q_p(2)\ (\mo\ p).
\endalign$$
So we also have
$$2q_p(2)-4\f{P_{p-(\f2p)}}p\eq2\sum_{k=1}^{\lfloor 3p/4\rfloor}\f{(-1)^{k-1}}k\ (\mo\ p).$$

(iii) Now suppose $p>3$. By Theorem 1.2 in the case
$m=2$, we have
$$\sum_{k=1}^{p-1}(-1)^k\f{\bi{2k}k}{k2^{k}}\eq\f{2^p-V_p(2)}p\
(\mo\ p).$$ Observe that the two roots of the equation
$x^2-2x-2=0$ are $1\pm\sqrt3$. Thus
$$\align V_p(2)=&(1+\sqrt3)^p+(1-\sqrt3)^p=2\sum_{k=0}^{(p-1)/2}\bi p{2k}(\sqrt3)^{2k}
\\=&2+\sum_{k=1}^{(p-1)/2}\f{2p}{2k}\bi{p-1}{2k-1}3^k
\\\eq&2-p\sum_{k=1}^{(p-1)/2}\f{3^k}{k}\ (\mo\ p^2).
\endalign$$
As observed by Eisenstein [2],
$$2q_p(2)=\f{2^p-2}p=\sum_{k=1}^{p-1}\f1p\bi
pk=\sum_{k=1}^{p-1}\f{\bi{p-1}{k-1}}k\eq\sum_{k=1}^{p-1}\f{(-1)^{k-1}}k\
(\mo\ p).$$ By a congruence of Z. W. Sun [15],
$$\sum_{k=1}^{(p-1)/2}\f{3^k}k\eq\sum_{0<k<p/6}\f{(-1)^k}k\eq\sum_{0<k<p/6}\f{(-1)^{p-k}}{p-k}\
(\mo\ p).$$ Thus
$$\align&\sum_{k=1}^{p-1}(-1)^k\f{\bi{2k}k}{k2^{k}}
\\\eq&\f{2^p-2}p-\f{V_p(2)-2}p\eq2q_p(2)+\sum_{k=1}^{(p-1)/2}\f{3^k}k
\\\eq&\sum_{k=1}^{p-1}\f{(-1)^{k-1}}k+\sum_{5p/6<k<p}\f{(-1)^k}k
=\sum_{0<k<5p/6}\f{(-1)^{k-1}}k\ (\mo\ p).
\endalign$$
In light of [15],
$$\sum_{k=1}^{(p-1)/2}\f{3^k}k\eq-q_p(2)-6\l(\f 2p\r)\f{S_{(p-(\f 3p))/2}}p\ (\mo\ p).$$
So we also have
$$\sum_{k=1}^{p-1}(-1)^k\f{\bi{2k}k}{k2^{k}}\eq
q_p(2)-6\l(\f 2p\r)\f{S_{(p-(\f 3p))/2}}p\ (\mo\ p).$$
Therefore (1.18) follows.

Let $u_n=u_n(2,-2)$ and $v_n=v_n(2,-2)$ for $n\in\N$. By induction,
$$v_n=2u_{n+1}-2u_n=2u_n+4u_{n-1}\quad\t{for}\ n=1,2,3,\ldots.$$
Thus
$$v_p=2\l(\f 3p\r)u_p+\l(3+\l(\f 3p\r)\r)u_{p-(\f 3p)}.$$
Clearly
$$\align 2\sqrt3u_{p-(\f 3p)}=&(1+\sqrt3)^{p-(\f 3p)}-(1-\sqrt3)^{p-(\f 3p)}
\\=&2^{(p-(\f 3p))/2}\l((2+\sqrt3)^{(p-(\f 3p))/2}-(2-\sqrt3)^{(p-(\f 3p))/2}\r)
\endalign$$
and hence
$$u_{p-(\f 3p)}=2^{(p-(\f 3p))/2}S_{(p-(\f 3p))/2}\eq\l(\f2p\r)2^{(1-(\f 3p))/2}S_{(p-(\f3p))/2}\ (\mo\ p^2).$$
Recall that
$$v_p=V_p(2)\eq2-p\sum_{k=1}^{(p-1)/2}\f{3^k}k\eq2+(2^{p-1}-1)+6\l(\f 2p\r)S_{(p-(\f3p))/2}\ (\mo\ p^2).$$
Therefore
$$\align&2\l(\f 3p\r)u_p-2=v_p-2-\l(3+\l(\f 3p\r)\r)u_{p-(\f 3p)}
\\\eq&2^{p-1}-1+6\l(\f 2p\r)S_{(p-(\f 3p))/2}-\l(3+\l(\f 3p\r)\r)\l(\f2p\r)2^{(1-(\f 3p))/2}S_{(p-(\f3p))/2}
\\\eq&2^{p-1}-1+2\l(\f 2p\r)S_{(p-(\f 3p))/2}\ (\mo\ p^2).
\endalign$$

Applying (1.11) with $m=-6$, we get
$$\align\f12\sum_{k=1}^{p-1}\f{\bi{2k}k}{k6^{k-1}}\eq&\f{V_p(-6)+6^p}p
\eq\f{V_p(-6)+6}p+6(q_p(2)+q_p(3))\ (\mo\ p).
\endalign$$
Observe that
$$\align V_p(-6)=&(-3+\sqrt3)^p+(-3-\sqrt3)^p
\\=&-\sqrt3^p\l((1+\sqrt3)^p-(1-\sqrt3)^p\r)
=-2\times 3^{(p+1)/2}u_p
\endalign$$
and hence
$$\align V_p(-6)+6\eq&-6\l(3^{(p-1)/2}-\l(\f 3p\r)\r)u_p-6\l(\f 3p\r)u_p+6
\\\eq&-6\l(3^{(p-1)/2}-\l(\f 3p\r)\r)\l(\f 3p\r)
\\&-3\l(2^{p-1}-1+2\l(\f 2p\r)S_{(p-(\f3p))/2}\r)
\\\eq&-3\l((3^{p-1}-1)+2^{p-1}-1+2\l(\f 2p\r)S_{(p-(\f3p))/2}\r)\ (\mo\ p^2).
\endalign$$
Therefore
$$\align \f12\sum_{k=1}^{p-1}\f{\bi{2k}k}{k6^{k-1}}\eq&
-3\(q_p(3)+q_p(2)+2\l(\f 2p\r)\f{S_{(p-(\f3p))/2}}p\)
\\&+6(q_p(2)+q_p(3))
\\\eq&3\(q_p(2)+q_p(3)-2\l(\f 2p\r)\f{S_{(p-(\f3p))/2}}p\)\ (\mo\ p).
\endalign$$
So (1.19) is valid.

 The proof of Corollary 1.3 is now complete. \qed

\heading{5. Proof of Theorem 1.3}\endheading

\medskip
\noindent{\it Proof of Theorem 1.3}. By an identity of T. B. Staver [12],
$$\sum_{k=1}^{n}\f{1}{k}\bi{2k}{k}=\f{2n+1}{3n^2}\bi{2n}{n}\sum_{k=1}^{n}\f1{\bi{n-1}{k-1}^{2}}
=\f{n+1}3\bi{2n+1}n\sum_{k=1}^n\f1{k^2\bi nk^2}$$ for all
$n=1,2,3,\ldots$. Taking $n=p^a-1$ in the identity, we get
$$\sum_{k=1}^{p^a-1}\f{1}{k}\bi{2k}{k}=\f{p^a}{3}\bi{2p^a-1}{p^a-1}
\sum_{k=1}^{p^a-1}\f1{k^2\bi{p^a-1}{k}^{2}}.\tag5.1$$

Recall that
$$\bi{2p^a-1}{p^a-1}\eq1+p[p=2]+p^2[p=3]\ (\mo\ p^3)$$
by [16, Lemma 2.2].  For  $k=1,\ldots,p^a-1$,  we set
$H_k=\sum_{0<j\ls k}1/j$ and note that
$$\align\f1{\bi{p^a-1}k^2}=&\prod_{0<j\ls k}\f1{(1-p^a/j)^2}
\\\eq&\prod_{0<j\ls k}\f{(1-p^{3a}/j^3)^2}{(1-p^a/j)^2}=\prod_{0<j\ls
k}\l(1+\f{p^a}j+\f{p^{2a}}{j^2}\r)^2
\\\eq&\prod_{0<j\ls k}\l(1+2\f{p^a}j+\f{p^{2a}}{j^2}+2\f{p^{2a}}{j^2}\r)
\ (\mo\ p^3)
\\\eq&\prod_{0<j\ls k}\l(1+2\f{p^a}j\r)\ (\mo\ p^{2+[p=3]})\endalign$$
Therefore (5.1) implies that
$$\align &p^{a-1}\sum_{k=1}^{p^a-1}\f{\bi{2k}k}k=\f p3\bi{2p^a-1}{p^a-1}\sum_{k=1}^{p^a-1}\f{p^{2(a-1)}}{k^2\bi{p^a-1}k^2}
\\\eq& \f p3(1+p[p=2]+p^2[p=3])\sum_{k=1}^{p^a-1}\f{p^{2(a-1)}}{k^2}\prod_{0<j\ls k}\l(1+2\f{p^a}j\r)\ (\mo\ p^3).
\endalign$$
So we have
$$p^{a-1}\sum_{k=1}^{p^a-1}\f{\bi{2k}k}k
\eq\l(\f p3+p^2[p\ls3]\r)\sum_{k=1}^{p^a-1}\f{p^{2(a-1)}}{k^2}\prod_{0<j\ls k}\(1+2\f{p^a}j\)\ (\mo\ p^3).\tag5.2$$

For $k=1,\ldots,p^a-1$, clearly
$$\align\prod_{0<j\ls k}\l(1+2\f{p^a}j\r)\eq&1+2p^aH_k+4p^{2a}\sum_{0<i<j\ls k}\f1{ij}
\\\eq&1+2p^aH_k+2p^{2a}\(H_k^2-\sum_{j=1}^k\f1{j^2}\)\ (\mo\ p^3).
\endalign$$
In the case $a\gs2$, if $1\ls k\ls p^a-1$ and $p^{a-2}\nmid k$ then $p^{2(a-1)}/k^2\eq0\ (\mo\ p^4)$.
When $a\gs 2$ and $k\in\{1,\ldots,p^2-1\}$, we have
$$\align\prod_{j=1}^{p^{a-2}k}\l(1+2\f{p^a}j\r)
\eq&1+2\sum_{j=1}^{p^{a-2}k}\f{p^a}j+2\(\sum_{j=1}^{p^{a-2}k}\f{p^a}j\)^2-2\sum_{j=1}^{p^{a-2}k}\f{p^{2a}}{j^2}
\\\eq&1+2\sum_{i=1}^k\f{p^a}{p^{a-2}i}+2\(\sum_{i=1}^k\f{p^a}{p^{a-2}i}\)^2-2\sum_{i=1}^k\f{p^{2a}}{(p^{a-2}i)^2}
\\\eq&1+2p^2H_k+2(p^2H_k)^2-2\sum_{i=1}^k\f{p^4}{i^2}\ (\mo\ p^3).
\endalign$$
Therefore, if $a\gs2$ then (5.2) implies that
$$\align p^{a-1}\sum_{k=1}^{p^a-1}\f{\bi{2k}k}k\eq&\l(\f p3+p^2[p\ls3]\r)\sum_{k=1}^{p^2-1}\f{p^{2(a-1)}}{(p^{a-2}k)^2}
\prod_{j=1}^{p^{a-2}k}\l(1+2\f{p^a}j\r)
\\\eq&\l(\f p3+p^2[p\ls3]\r)\sum_{k=1}^{p^2-1}\f{p^2}{k^2}\prod_{j=1}^k\l(1+2\f{p^2}j\r)\ (\mo\ p^3).
\endalign$$
In the case $p=3$, this yields (1.20) for $a\gs2$. (1.20) in the case $p=3$ and $a=1$ can be verified directly.

 Below we assume that $p\not=3$.
For $k=1,\ldots,p^a-1$, if $p^{a-1}\nmid k$ then $p^{2(a-1)}/k^2\eq0
(\mo\ p^2)$. Also,
$$p^aH_{p^{a-1}k}=\sum_{j=1}^{p^{a-1}k}\f{p^a}j\eq\sum_{i=1}^k\f{p^a}{p^{a-1}i}=pH_k\ (\mo\ p^2)$$
for every $k=1,\ldots,p-1$.
Thus (5.2) implies that
$$\align p^{a-1}\sum_{k=1}^{p^a-1}\f{\bi{2k}k}k\eq&
\l(\f p3+p^2[p=2]\r)\sum_{k=1}^{p-1}\f{p^{2(a-1)}}{(p^{a-1}k)^2}(1+2p^aH_{p^{a-1}k})
\\\eq& \l(\f p3+p^2[p=2]\r)\sum_{k=1}^{p-1}\f{1+2pH_k}{k^2}\ (\mo\ p^3).
\endalign$$
This yields (1.20) in the case $p=2$.

 Now we handle the remaining case $p>3$. By the above, it suffices to show that
 $$\sum_{k=1}^{p-1}\f{1+2pH_k}{k^2}\eq\f 83pB_{p-3}\ (\mo\ p^2).\tag5.3$$
Let $n\in\N$. It is well known that
$$\sum_{j=0}^{k-1}j^n=\f1{n+1}\sum_{i=0}^n\bi{n+1}iB_ik^{n+1-i}\quad\ \t{for}\ k\in\Z^+,$$
and that
$$\sum_{k=1}^{p-1}k^n\eq pB_n\eq0\ (\mo\ p)\ \ \  \t{if}\ n\not\eq0\ (\mo\ p-1).$$
(See, e.g., [6, p.\,235].) Therefore
$$\align
\sum_{k=1}^{p-1}\f1{k^2}\sum_{j=0}^kj^{p-2}
=&\sum_{k=1}^{p-1}\(k^{p-4}+\f1{k^2(p-1)}\sum_{i=0}^{p-2}\bi{p-1}iB_ik^{p-1-i}\)
\\=&\sum_{k=1}^{p-1}k^{p-4}+\f1{p-1}\sum_{i=0}^{p-2}\bi{p-1}iB_i\sum_{k=1}^{p-1}k^{p-3-i}
\\\eq&\bi{p-1}{p-3}B_{p-3}+\f{B_{p-2}}2\sum_{k=1}^{p-1}\l(\f1k+\f1{p-k}\r)\ (\mo\ p)
\endalign$$
and hence
$$\sum_{k=1}^{p-1}\f{H_k}{k^2}\eq B_{p-3}\ (\mo\ p).\tag5.4$$
By a result of J. W. L. Glaisher [3, 4],
$$\bi{2p-1}{p-1}\eq 1-p^2\sum_{k=1}^{p-1}\f1{k^2}\eq1-\f23p^3B_{p-3}\ (\mo\ p^4)$$
and thus
$$\sum_{k=1}^{p-1}\f1{k^2}\eq\f 23pB_{p-3}\ (\mo\ p^2).\tag5.5$$
Note that (5.3) follows from (5.4) and (5.5).
We are done. \qed

\enddocument